\documentclass[11pt,leqno]{amsart2000}
\newtheorem{theo}{Theorem}[section]
\newtheorem{lm}[theo]{Lemma}
\newtheorem{cor}[theo]{Corollary}
\newtheorem{pro}[theo]{Proposition}
\theoremstyle{definition}
\newtheorem{exa}[theo]{Example}
\newtheorem{exas}[theo]{Examples}
\newtheorem{defi}[theo]{Definition}

\newtheorem{notas}[theo]{Notations}
\newtheorem{rem}[theo]{Remark}

\newtheorem{fact}[theo]{Fact}

\newtheorem{nist}[theo]{}
\def\a{\alpha}
\def\b{\beta}
\def\g{\gamma}
\def\d{\delta}
\def\ep{\varepsilon}
\def\La{\Leftarrow}
\def\Ra{\Rightarrow}
\def\lra{\to}
\def\sbe{\subseteq}
\def\spe{\supseteq}
\def\stm{\setminus}
\def\ems{\emptyset}
\def\ovl{\overline}
\def\ap{^\prime}
\def\card #1{\vert #1 \vert}
\def\AA{\mathcal{A}}
\def\BB{\mathcal{B}}
\def\CC{\mathcal{C}}
\def\DD{\mathcal{D}}
\def\FF{\mathcal{F}}
\def\HH{\mathcal{H}}
\def\MM{\mathcal{M}}
\def\OO{\mathcal{O}}
\def\PP{\mathcal{P}}
\def\TT{\mathcal{T}}
\def\UU{\mathcal{U}}
\def\Z{{\mathbb Z}}
\def\RRRR{{\mathbb R}}
\def\NNNN{{\mathbb N}}
\def\Ap{A^+_\mathcal{M}}

\begin{document}
\setlength{\unitlength}{0.01in}
\linethickness{0.01in}
\begin{center}
\begin{picture}(474,66)(0,0)
\multiput(0,66)(1,0){40}{\line(0,-1){24}}
\multiput(43,65)(1,-1){24}{\line(0,-1){40}}
\multiput(1,39)(1,-1){40}{\line(1,0){24}}
\multiput(70,2)(1,1){24}{\line(0,1){40}}
\multiput(72,0)(1,1){24}{\line(1,0){40}}
\multiput(97,66)(1,0){40}{\line(0,-1){40}}
\put(143,66){\makebox(0,0)[tl]{\footnotesize Proceedings of the Ninth Prague Topological Symposium}}
\put(143,50){\makebox(0,0)[tl]{\footnotesize Contributed papers from the symposium held in}}
\put(143,34){\makebox(0,0)[tl]{\footnotesize Prague, Czech Republic, August 19--25, 2001}}
\end{picture}
\end{center}
\vspace{0.25in}
\setcounter{page}{51}
\title{On Tychonoff-type hypertopologies}
\author{Georgi Dimov}
\address{Department of Mathematics\\
University of Sofia\\
Blvd. J. Bourchier 5\\
1126 Sofia, Bulgaria}
\email{gdimov@fmi.uni-sofia.bg}
\author{Franco Obersnel}
\author{Gino Tironi}
\address{Department of Mathematical Sciences and Computer Science\\
University of Trieste\\
Via A. Valerio 12/1\\
34127 Trieste, Italy}
\email{obersnel@mathsun1.univ.trieste.it}
\email{tironi@univ.trieste.it}
\begin{abstract}
In 1975, M. M. Choban \cite{C} introduced a new topology on the set of all
closed subsets of a topological space, similar to the {\em Tychonoff
topology} but weaker than it. In 1998, G. Dimov and D. Vakarelov
\cite{DV} used a generalized version of this new topology, calling it {\em
Tychonoff-type topology}. 
The present paper is devoted to a detailed study of Tychonoff-type 
topologies on an arbitrary family $\MM$ of subsets of a set $X$. 
When $\MM$ contains all singletons, a description of all Tychonoff-type 
topologies $\OO$ on $\MM$ is given.
The continuous maps of a special form between spaces of the type
$(\MM,\OO)$ are described in an isomorphism theorem.
The problem of {\em commutability between hyperspaces and subspaces with
respect to a Tychonoff-type topology} is investigated as well.
Some topological properties of the hyperspaces $(\MM,\OO)$ with
Tychonoff-type topologies $\OO$ are briefly discussed.
\end{abstract}
\thanks{The first author was partially supported by a Fellowship for
Mathematics of the NATO-CNR Outreach Fellowships Programme 1999, Bando
219.32/16.07.1999.}
\thanks{The second and the third authors were supported by the National
Group ``Analisi reale'' of the Italian Ministry of the University and
Scientific Research at the University of Trieste.}
\thanks{Georgi Dimov, Franco Obersnel and Gino Tironi,
{\em On Tychonoff-type hypertopologies},
Proceedings of the Ninth Prague Topological Symposium, (Prague, 2001),
pp.~51--70, Topology Atlas, Toronto, 2002}
\maketitle

\section{Introduction}

In 1975, M. M. Choban \cite{C} introduced a new topology on the set of all 
closed subsets of a topological space for obtaining a generalization of
the famous Kolmogoroff Theorem on operations on sets.
This new topology is similar to the {\em Tychonoff topology} (known also
as {\em upper Vietoris topology}, or {\em upper semi-finite topology}
(\cite{M}), or {\em kappa-topology}) but is weaker than it. 
In 1998, G. Dimov and D. Vakarelov \cite{DV} used a generalized version of
this new topology for proving an isomorphism theorem for the category of
all Tarski consequence systems. 
This generalized version was called {\em Tychonoff-type topology}.

The present paper is devoted to a detailed study of Tychonoff-type
topologies on an arbitrary family $\MM$ of subsets of a set $X$. 
When $\MM$ is a {\em natural family}, i.e.\ it contains all singletons, a
description of all Tychonoff-type topologies $\OO$ on $\MM$ is given (see
Proposition~\ref{char Tych}).
For doing this, the notion of {\em T-space} is introduced. 
The natural morphisms for T-spaces are not enough to describe all
continuous maps between spaces of the type $(\MM,\OO)$, where $\MM$ is a
natural family and $\OO$ is a Tychonoff-type topology on it; we obtain a
characterization of those continuous maps which correspond to the
morphisms between T-spaces.
This is done by defining suitable categories and by proving that these
categories are isomorphic (see Theorem~\ref{isomorphism}).
In such a way we extend to any natural family $\MM$ on $X$ the 
corresponding result obtained in \cite{DV} for the family $\FF in(X)$
of all finite subsets of $X$. 
We investigate also the problem of {\em commutability between hyperspaces
and subspaces with respect to a Tychonoff-type topology}, i.e.\ when the
hyperspace of any subspace $A$ of a topological space $Y$ is canonically
representable as a subspace of the hyperspace of $Y$.
Such investigations were done previously by H.-J. Schmidt \cite{S} for the
lower Vietoris topology, by G. Dimov \cite{D1, D2} for the Tychonoff 
topology and for the Vietoris topology, and by B. Karaivanov \cite{K} for
other hypertopologies. 
We study also such a problem for a fixed subspace $A$ of $Y$. 
Some results of \cite{D1, D2, Se} are generalized.
Finally, we study briefly some topological properties (separation axioms, 
compactness, weight, density, isolated points, $P_\infty$) of the 
hyperspaces $(\MM,\OO)$ with Tychonoff-type topologies $\OO$.
Some results of \cite{F, DV} are generalized.

Let us fix the notations.

\begin{notas}\label{initial notation}
We denote by $\omega$ the set of all {\em positive} natural numbers, by
$\RRRR$ --- the real line, and by $\Z$ --- the set of all integers.
We put $\NNNN=\omega\cup\{0\}$.

Let $X$ be a set. 
We denote by $\PP(X)$ the set of all subsets of $X$.
Let $\MM, \AA \subseteq\PP(X)$ and $A\subseteq X$.
We will use the following notations:
\begin{itemize}
\item
$A_\MM^+:=\{M\in\MM : M\subseteq A\};$
\item
$\AA^+_\MM:=\{A^+_\MM : A\in\AA\};$
\item
$\FF in (X):=\{M\subseteq X : |M|<\aleph_0\}$;
\item
$\FF in_n (X):=\{M\subseteq X : |M|\le n\}$, where $n\in\omega$.
\end{itemize}
We will denote by $\AA^\cap$ (respectively by $\AA^\cup$) the closure 
under finite intersections (unions) of the family $\AA$.
In other words,
\begin{itemize}
\item
$\AA^\cap:=\{\bigcap_{i=1}^k A_i : k\in \omega, A_i\in\AA\}$ and
\item
$\AA^\cup:=\{\bigcup_{i=1}^k A_i : k\in \omega, A_i\in\AA\}.$
\end{itemize}

Let $(X,\TT)$ be a topological space. 
We put 
\begin{itemize}
\item
$\CC L(X):=\{ M\subseteq X: M$ is closed in $X,\ M\not=\emptyset\}$
and
\item
$\CC omp(X):=\{ M\subseteq X: M { \rm \ is \ compact}\}$.
\end{itemize}
The closure of a subset $A$ of $X$ in $(X,\TT)$ will be denoted by $cl_X
A$ or $\ovl{A}^X$; as usual, for $U\sbe A\sbe X$, we put 
\begin{itemize}
\item
$Ex_{A,X}U:=X\stm cl_X(A\stm U)$.
\end{itemize}
By a {\it base} of $(X,\TT)$ we will always mean an open base.
The weight (resp., the density) of $(X,\TT)$ will be denoted by $w(X,\TT)$
(resp., $d(X,\TT)$).

If $\CC$ denotes a category, we write $X\in |\CC|$ if $X$ is an object of
the category $\CC$.

For all undefined here notions and notations, see \cite{E} and \cite{J}.
\end{notas}

\section{Hypertopologies of Tychonoff-type}

\begin{fact}\label {intersection}
Let $X$ be a set and $\MM,\AA \subseteq\PP(X)$. 
Then:
\begin{itemize}
\item[(a)] 
$\bigcap \AA^+_\MM = \left(\bigcap \AA\right)^+_\MM;$
\item[(b)] 
$A\subseteq B$ implies that $\Ap\subseteq B^+_\MM$ for all 
$A$, $B \subseteq X$.
\end{itemize}
\end{fact}

\begin{defi}\label{Tych top}
Let $(X,\TT)$ be a topological space and let $\MM\subseteq \PP(X)$. 
The topology $\OO_\TT$ on $\MM$, having as a base the family $\TT^+_\MM$,
will be called {\it Tychonoff topology on $\MM$ generated by 
$(X,\TT)$}. When $\MM=\CC L(X)$, then $\OO_\TT$ is just the classical 
Tychonoff topology on $\CC L(X)$.

Let $X$ be a set and $\MM\subseteq \PP(X)$. 
A topology $\OO$ on $\MM$ is called a {\it Tychonoff topology on} $\MM$ if
there exists a topology $\TT$ on $X$ such that $\TT^+_\MM$ is a base of
$\OO$.
\end{defi}

\begin{defi}\label{Tych type top}
Let $X$ be a set and $\MM\subseteq\PP (X)$. 
A topology $\OO$ on the set $\MM$ is called a {\it topology of
Tychonoff-type on $\MM$} if the family $\OO\cap\PP(X)^+_\MM$ is a base for 
$\OO$.
\end{defi}

Clearly, a Tychonoff topology on $\MM$ is always a topology of 
Tychonoff-type on $\MM$, but not viceversa (see Example~\ref{example}).

\begin{fact}\label{BO}
Let $X$ be a set, $\MM\subseteq\PP (X)$ and $\OO$ be a topology of
Tychonoff-type on $\MM$.
Then the family $\BB_\OO:=\{A\subseteq X : \Ap\in\OO\}$ is closed under
finite intersections, $X\in \BB_\OO$, and, hence, $\BB_\OO$ is a base for
a topology $\TT_\OO$ on $X$.
The family $\left(\BB_\OO\right)^+_\MM$ is a base of $\OO$.
\end{fact}

\begin{defi}\label{induced}
Let $X$ be a set, $\MM\subseteq\PP (X)$ and $\OO$ be a topology of 
Tychonoff-type on $\MM$. 
We will say that the topology $\TT_\OO$ on $X$, introduced in 
Fact~\ref{BO}, is {\it induced by the topological space $(\MM,\OO)$.}
\end{defi}

\begin{pro}\label{char top Tych type}
Let $X$ be a set and $\MM\subseteq\PP(X)$. 
A topology $\OO$ on $\MM$ is a topology of Tychonoff-type if and only if
there exists a topology $\TT$ on $X$ and a base $\BB$ for $\TT$ (which
contains $X$ and is closed under finite intersections) such that
$\BB^+_\MM$ is a base for $\OO$.
\end{pro}

\begin{proof}
Suppose $\OO$ is a topology of Tychonoff-type on $\MM$. 
Then the topology $\TT_\OO$ induced by the topological space $(\MM,\OO)$
(see Fact~\ref{BO} and Definition~\ref{induced}) and the base $\BB_\OO$
have the required property.

Conversely, suppose $\TT$ and $\BB$ are given as in the statement. 
Then $\BB^+_\MM$ is a base for $\OO$, and therefore also 
$\OO\cap\PP(X)^+_\MM$ is a base for $\OO$.
\end{proof}

\begin{defi}\label{generates}
Let $X$ be a set and $\MM, \BB \subseteq\PP (X)$.
When $\BB_\MM^+$ is a base for a topology $\OO_\BB$ on $\MM$, we will say
that $\BB$ {\em generates a topology} on $\MM$.
(Obviously, the topology $\OO_\BB$ is of Tychonoff-type.
\end{defi}

\begin{pro}\label{char base}
Let $X$ be a set, $\MM\subseteq\PP(X)$ and $\BB\subseteq \PP(X)$.
The family $\BB$ generates a topology $\OO_\BB$ on $\MM$ if and only if
the family $\BB$ satisfies the following conditions:
\begin{itemize}
\item[(MB1)] 
For any $M\in\MM$ there exists a $U\in \BB$ such that $M\subseteq U$;
\item[(MB2)] 
For any $U_1$, $U_2 \in \BB$ and any $M\in\MM$ with 
$M\subseteq U_1\cap U_2$ there exists a $U_3\in\BB$ such that
$M\subseteq U_3\subseteq U_1\cap U_2$.
\end{itemize}
\end{pro}

\begin{proof}
It follows from Proposition 1.2.1 \cite{E}.
\end{proof}

\begin{cor}\label{cord}
Let $X$ be a set and $\MM, \BB\sbe \PP(X)$. 
If $\BB=\BB^\cap$ and $X\in \BB$, then $\BB$ generates a Tychonoff-type
topology on $\MM$.
\end{cor}

\begin{defi}\label{natural family}
Let $X$ be a set and $\MM\subseteq\PP(X)$. 
We say that $\MM$ is a {\it natural family in $X$} if $\{x\}\in\MM$ for
all $x\in X$.
\end{defi}

\begin{cor}\label{base on X}
Let $X$ be a set and $\MM$ be a natural family in $X$. 
If $\BB\subseteq\PP(X)$ generates a topology on $\MM$ (see Definition
\ref{generates}), then $\BB$ is a base for a topology on $X$.
\end{cor}

\begin{proof}
By Proposition~\ref{char base}, $\BB$ satisfies the conditions (MB1) and
(MB2). 
Since $\MM$ is natural, this clearly implies that $\BB$ satisfies the
hypotesis of Proposition 1.2.1 \cite{E}.
So $\BB$ is a base for a topology on $X$.
\end{proof}

\begin{rem}
Trivial examples show that there exist sets $X$ and (non-natural) families
$\MM$,$\BB\subseteq\PP (X)$ such that $\BB_\MM^+$ is a base for a topology
on $\MM$ but
\begin{itemize}
\item[(a)] 
$\bigcup\BB\not= X$, so that $\BB$ cannot serve even as subbase of a
topology on $X$ (take $X=\{0,1\}$, $\MM=\BB=\{\{0\}\}$);
\item[(b)]
$\BB$ is not a base of a topology of $X$, although $\bigcup\BB= X$
(take $X=\{0,1,2\}$, $\MM=\BB=\{\{0,1\},\{0,2\}\}$).
\end{itemize}
The example of (b) shows also that if we substitute in \ref{base on X} 
naturality of $\MM$ with the condition ``$\bigcup\MM=X$'' then we cannot
prove that $\BB$ is a base of a topology on $X$; however, it is easy to
show that the condition ``$\bigcup\MM=X$'' implies that $\bigcup\BB=X$,
i.e.\ $\BB$ can serve as a subbase of a topology on $X$.

Of course, as it follows from Fact~\ref{intersection}, if $\BB^+_\MM$ is a
base of a topology $\OO$ on $\MM$, then $\tilde \BB=\BB^\cap \cup\{X\}$ is
a base for a topology on $X$ and $\tilde\BB_\MM^+$ is a base of $\OO$.
\end{rem}

\begin{cor}\label{Fin1}
Let $(X,\TT)$ be a topological space and let $\BB\subseteq\TT$ be a base of
$(X,\TT)$, closed under finite unions. 
Then $\BB$ generates a topology of Tychonoff-type on $\FF in(X)$ and 
$\CC omp(X)$.
\end{cor}

\begin{proof}
It follows easily from Proposition~\ref{char base}.
\end{proof}

\begin{pro}\label{O1=O2}
Let $X$ be a set, $\MM$,$\BB_1$,$\BB_2 \subseteq\PP(X)$, and suppose that
$\BB_1$ and $\BB_2$ generate, respectively, some topologies (of 
Tychonoff-type) $\OO_{\BB_1}$ and $\OO_{\BB_2}$ on $\MM$.
Then $\OO_{\BB_1}=\OO_{\BB_2}$ if and only if the following conditions are
satisfied:
\begin{itemize}
\item[(CO1)] 
For any $M\in\MM$ and any $U_1\in \BB_1$ such that $M\subseteq U_1$ there
exists $U_2\in \BB_2$ with $M\subseteq U_2\subseteq U_1$;
\item[(CO2)] 
For any $M\in\MM$ and any $U_2\in \BB_2$ such that $M\subseteq U_2$ there
exists $U_1\in \BB_1$ with $M\subseteq U_1\subseteq U_2$.
\end{itemize}
\end{pro}

\begin{proof}
It follows from 1.2.B \cite{E}.
\end{proof}

\begin{cor}\label {Fin2}
Let $(X,\TT)$ be a topological space and $\BB_1$, $\BB_2$ be bases of 
$(X,\TT)$, closed under finite unions. 
Then they generate (see Corollary~\ref{Fin1}) equal topologies on 
$\FF in(X)$ and $\CC omp(X)$. 
In particular, every topology of Tychonoff-type on $\FF in(X)$ or on 
$\CC omp(X)$, generated by a base of $(X,\TT)$ which is closed under
finite unions, coincides with the Tychonoff topology generated by
$(X,\TT)$ on the corresponding set.
\end{cor}

\begin{proof}
Check that conditions (CO1) and (CO2) of Proposition~\ref{O1=O2} are
satisfied.
\end{proof}

\begin{cor}\label{Mbase}
Let $(X,\TT)$ be a topological space, $\MM\subseteq\PP(X)$, 
$\BB\subseteq\TT$, and suppose that $\BB$ generates a topology of
Tychonoff-type $\OO$ on $\MM$.
Then $\OO$ is the Tychonoff topology on $\MM$ generated by $(X,\TT)$ if
and only if for all $M\in\MM$ and for all $V\in\TT$ such that 
$M\subseteq V$, there exists $U\in\BB$ with $M\subseteq U\subseteq V$. 
In this case we will say that $\BB$ is an {\em $\MM$-base} for $(X,\TT)$.
Clearly, if $\MM$ is a natural family, then every $\MM$-base of $(X,\TT)$
is also a base of $(X,\TT)$.
\end{cor}

\begin{proof}
Put $\BB_1 :=\TT$ and $\BB_2 :=\BB$. 
Then condition (CO2) of Proposition~\ref{O1=O2} is trivially satisfied. 
The condition required in the statement is exactly condition (CO1).
\end{proof}

\begin{defi}\label{Mcover}
Let $X$ be a set, $\MM\subseteq\PP(X)$ and $A\subseteq X$. 
A family $\UU\subseteq\PP(X)$ will be called an {\it $\MM$-cover} of $A$
if $A=\bigcup \UU$ and for all $M\in\MM$ with $M\subseteq A$ there exists
some $U\in\UU$ such that $M\subseteq U$.
\end{defi}

\begin{pro}\label{union2}
Let $X$ be a set and $\MM, \AA\subseteq \PP(X)$. Then the following
conditions are equivalent:
\begin{itemize}
\item[(U1)]
For all $U\in \AA$ and for all $x\in U$, there exists an $M\in\MM$ with
$x\in M\subseteq U$.
\item[(U2)]
For any $U\in\AA\cup\MM$ and for any subfamily 
$\{U_\delta : \delta\in \Delta\}$ of $\AA\cup\MM$, the equality 
$U^+_\MM = \bigcup_{\delta\in\Delta} \left( U_\delta\right)^+_\MM$ 
holds if and only if the family $\{U_\delta\}_{\delta\in \Delta}$ is an
$\MM$-cover of $U$.
\end{itemize}
\end{pro}

\begin{proof}
Observe that, trivially, in condition (U1) we can replace the requirement
`for all $U\in\AA$'' with ``for all $U\in\AA\cup\MM$''.

(U1)$\Rightarrow$(U2).
Let $U^+_\MM = \bigcup_{\delta\in\Delta} \left( U_\delta\right)^+_\MM$, with
$U$,$U_\delta\in\AA\cup\MM$ for all $\d\in\Delta$.
We will prove first that $\bigcup_{\delta\in\Delta}U_\delta = U$.

Let $x\in\bigcup_{\delta\in\Delta}U_\delta$. 
Then there exists a $\delta\in\Delta$ such that $x\in U_\delta$. 
By assumption, there exists an $M\in\MM$ with 
$x\in M\subseteq U_\delta$. 
Hence $M\in \left( U_\delta\right)^+_\MM$. 
Since $U^+_\MM = \bigcup_{\delta\in\Delta} \left( U_\delta\right)^+_\MM$,
we obtain that $M\subseteq U$. 
Thus $x\in U$. 
Therefore, $\bigcup_{\delta\in\Delta}U_\delta \sbe U$.

Conversely, let $x\in U$. 
By assumption, there exists an $M\in\MM$ such that 
$x\in M\subseteq U$. 
Hence 
$M\in U^+_\MM = \bigcup_{\delta\in\Delta} \left( U_\delta\right)^+_\MM$.
Therefore there exists a $\delta\in\Delta$ such that 
$M\in \left( U_\delta\right)^+_\MM$, i.e.\ $M\subseteq U_\delta$ and 
$x\in U_\delta \subseteq \bigcup_{\delta\in\Delta} U_\delta$.

We have verified that $\bigcup_{\delta\in\Delta}U_\delta = U$.

Suppose $M\in\MM$ and $M\subseteq U$. 
Then $M\in U^+_\MM $ and therefore there exists some $\gamma\in\Delta$
with $M\in \left( U_\gamma\right)^+_\MM$.
Hence $M\sbe U_\g$.

This shows that the family $\{U_\delta\}_{\delta\in \Delta}$ is an
$\MM$-cover of $U$.

The other implication can be easily proved.
(Let's remark that condition (U1) is not used in the proof of this last
implication.)

(U2)$\Rightarrow$(U1).
Suppose $U\in\AA$ and $x\in U$.
Clearly, we have 
$$U^+_\MM=\bigcup\{ M^+_\MM: {M\in\MM,\ M\subseteq U}\}.$$
Then, by assumption, the family $\{M\in\MM:M\subseteq U\}$ is an
$\MM$-cover of $U$.
Therefore $U=\bigcup\{ M: {M\in\MM,\ M\subseteq U}\}$. 
Hence there exists an $M\in\MM$ with $x\in M\subseteq U$.
\end{proof}

\begin{pro}\label{union}
Let $X$ be a set and $\MM\sbe\PP(X)$. Then the following
conditions are equivalent:
\begin{itemize}
\item[(a)]
$\MM$ is a natural family;
\item[(b)]
For any $U\sbe X$ and for any subfamily $\{U_\delta : \delta\in \Delta\}$
of $\PP(X)$, the equality 
$U^+_\MM = \bigcup_{\delta\in\Delta} \left(U_\delta\right)^+_\MM$ 
holds if and only if the family $\{U_\delta\}_{\delta\in \Delta}$ is an
$\MM$-cover of $U$.
\end{itemize}
\end{pro}

\begin{proof}
Put $\AA=\PP(X)$ in Proposition~\ref{union2}.
\end{proof}

\begin{pro}\label{uniqueness}
Let $X$ be a set, $\MM\subseteq\PP(X)$, $\OO$ be a Tychonoff topology on
$\MM$ generated by a topology $\TT$ on $X$ and
$\MM$ be a network in the sense of Arhangel$'$ski\u{\i}\ for $\TT_\OO$.
Then
$\TT=\BB_\OO$ and
$\OO$
is generated by a unique topology on $X$, namely by $\TT_\OO$
(see Fact~\ref{BO} for the notation $\BB_\OO$ and $\TT_\OO$).
\end{pro}

\begin{proof}
We only need to show that $\BB_\OO\subseteq\TT$. Assume $A\in \BB_O$. Then
$A^+_\MM\in\OO$. Since
$\TT$ generates $\OO$, we have
$A^+_\MM=\bigcup_{\delta\in\Delta}\left(U_\delta\right)^+_\MM$, where
$U_\delta\in\TT$ for all $\delta\in\Delta$.
Clearly $\ U_\delta\in\BB_O$ for all
$\delta\in\Delta$. By Proposition~\ref{union2}, we obtain that
$A=\bigcup_{\delta\in\Delta} U_\delta$ and
therefore $A\in\TT$.
\end{proof}

\begin{rem}\label{rema}
Trivial examples show that there exist sets $X$, families
$\MM\sbe \PP(X)$ and Tychonoff topologies on $\MM$ which are generated
by more than one topology on $X$.
\end{rem}

\begin{cor}\label{char Tych top}
Let $X$ be a set, $\MM\subseteq\PP(X)$, $\OO$ be a
topology of Tychonoff-type on $\MM$ and
$\MM$ be a network in the sense of Arhangel$'$ski\u{\i}\ for $\TT_\OO$.
Then $\OO$ is a Tychonoff topology on $\MM$ if and only if $\BB_\OO=\TT_\OO$
(see
Fact~\ref{BO}
for the notations $\BB_\OO$ and $\TT_\OO$).
\end{cor}

\begin{proof}
Suppose $\BB_\OO=\TT_\OO$. Then the topology $\OO$ is generated by the
topology
$\TT_\OO$ on $X$ and hence, by
definition, $\OO$ is a Tychonoff topology on $\MM$.

Suppose $\OO$ is a Tychonoff topology on $\MM$.
Then $\OO$ is generated by some topology $\TT$ on $X$.
By Proposition~\ref{uniqueness}, we get $\TT=\BB_\OO$. Hence
$\TT_\OO=\BB_\OO$.
\end{proof}

\begin{cor}\label{unique}
Let $X$ be a set, $\MM$ be a natural family in $X$ and $\OO$ be a topology
of Tychonoff-type on
$\MM$. Then $\OO$ is a Tychonoff topology on $\MM$ if and only if
$\BB_\OO=\TT_\OO$.
\end{cor}

\begin{proof}
A natural family $\MM$ satisfies the hypothesis of 
Corollary~\ref{char Tych top}.
\end{proof}

\begin{pro}\label {char BOB}
Let $(X,\TT)$ be a topological space, $\MM$ be a natural family in $X$,
$\BB\subseteq\TT$
and suppose that $\BB$ generates a topology $\OO_\BB$ on $\MM$.
Then
$$\BB_{\OO_\BB}=\{A\subseteq X : A \mbox{ \ is\ } \MM\mbox{-covered \ by\ a
\
subfamily \ of\ } \BB\},$$
$\BB_{\OO_\BB}\subseteq \TT$ and $\BB^\cap\subseteq \BB_{\OO_\BB}$ (see
Fact~\ref{BO} for the notation
$\BB_{\OO_\BB}$).
\end{pro}

\begin{proof}
It follows from Proposition~\ref{union} and Fact~\ref{BO}.
\end{proof}

\begin{pro}\label{TOB=T}
Let $(X,\TT)$ be a topological space, $\MM$ be a natural family in $X$,
$\BB$ be a base for
$(X,\TT)$
and suppose that $\BB$ generates a topology $\OO_\BB$ on $\MM$.
Let
$\TT_{\OO_\BB}$ be the topology on $X$ induced by $(\MM,\OO_\BB)$.
Then $\TT_{\OO_\BB}=\TT$.
\end{pro}

\begin{proof}
We have, by Proposition~\ref{char BOB},
that $\BB\subseteq \BB_{\OO_\BB}\subseteq \TT_{\OO_\BB}$.
Thus $\TT\subseteq\TT_{\OO_\BB}$.
As it is shown in Proposition~\ref{char BOB},
$\BB_{\OO_\BB}\subseteq \TT$ and
hence
$\TT_{\OO_\BB}\subseteq \TT$. So, $\TT=\TT_{\OO_\BB}$.
\end{proof}

\begin{exa}\label{ex1}
Let us show that
in Proposition~\ref{TOB=T}
the requirement
``$\MM$ is a natural family'' is essential.

Let $X=(0,1)\subset\RRRR$ be the open unit interval with the usual topology,
$\MM=\{[a,b] : 0<a<b<1\}$,
$\BB=\{(a,b): 0<a<b<1\}$. Then the family
$\BB$ satisfies conditions (MB1) and (MB2). Consider the set
$A=(\frac{1}{2},\frac{3}{4})\cup\{\frac{1}{4}\}$. We have $A^+_\MM =
(\frac{1}{2},\frac{3}{4})^+_\MM\in\BB^+_\MM$ and therefore
$A\in \BB_{\OO_\BB}$ even though $A$ is not open in $X$.
\end{exa}

\begin{defi}\label{Mclosed}
Let $X$ be a set and $\MM, \UU\sbe \PP(X)$. We will
say that $\UU$ is an {\it $\MM$-closed family} if for all $A\subseteq X$
such that
$A$ is $\MM$-covered by some subfamily of $\UU$, we have that $A\in\UU$.
\end{defi}

\begin{pro}\label{M and M^prime}
Let $X$ be a set and $\MM$, $\MM^\prime$, $\UU$ $\subseteq\PP(X)$,
$\MM\subseteq\MM^\prime$. Suppose that $\UU$ is
an $\MM$-closed family. Then $\UU$ is an $\MM^\prime$-closed family too.
\end{pro}

\begin{proof}
Suppose $A\subseteq X$ is $\MM^\prime$-covered by some subfamily
$\UU^\prime$
of $\UU$. Since
$\MM\subseteq\MM^\prime$, the set $A$ is also
$\MM$-covered by $\UU^\prime$. By the hypothesis, $\UU$ is an $\MM$-closed
family. Hence $A\in\UU$.
\end{proof}

\begin{pro}\label{BO Mclosed}
Let $X$ be a set, $\MM$ be a natural family in $X$ and $\OO$ be a topology
of Tychonoff-type on $\MM$. 
Then $\BB_\OO$ is an $\MM$-closed family (see Fact~\ref{BO} for the
notation $\BB_\OO$).
\end{pro}

\begin{proof}
It follows from Proposition ~\ref{union}.
\end{proof}

\begin{pro}\label{BOB=B}
Let $(X,\TT)$ be a topological space, $\MM$ be a natural family in $X$ and
$\BB\subseteq \TT$ be an $\MM$-closed base of $(X,\TT)$.
Suppose that $\BB$ generates a topology $\OO_\BB$ on $\MM$. 
Then $\BB_{\OO_\BB}=\BB$, $X\in\BB$ and $\BB^\cap =\BB$
(see Fact~\ref{BO} for the notation $\BB_{\OO_\BB}$).
\end{pro}

\begin{proof}
Obviously, $\BB_{\OO_\BB}\supseteq\BB$.
Let us show that $\BB_{\OO_\BB}\subseteq\BB$.

Let $A\in \BB_{\OO_\BB}$.
Then, by Proposition~\ref{char BOB}, $A$ is $\MM$-covered by some
subfamily of $\BB$ and, since $\BB$ is $\MM$-closed, we conclude that
$A\in\BB$. 
So, $\BB_{\OO_\BB}=\BB$. 
Now Fact~\ref{BO} implies that $X\in\BB$ and $\BB^\cap =\BB$.
\end{proof}

\begin{defi}\label{T_M-space}
Let $X$ be a set and $\MM, \BB\subseteq\PP(X)$.
The ordered triple $(X,\BB,\MM)$ will be called a {\em $T$-space} if
$\BB$ is an
$\MM$-closed family, $X\in\BB$
and $\BB^\cap=\BB$.
\end{defi}

Note that if $(X,\BB,\MM)$ is a $T$-space, then $\BB$ is a base for a 
topology on $X$.

\begin{pro}\label{char Tych}
Let $X$ be a set and $\MM$ be a natural family in $X$.
Let $TTT(X,\MM)$ be the set of all topologies of Tychonoff-type on $\MM$. 
Denote by $T\!\!{\rm-}\!\!Sp(X,\MM)$
the set of all $T$-spaces of the form $(X,\BB,\MM)$.
Then there is a bijective correspondence between the sets $TTT(X,\MM)$ and
$T\!\!{\rm -}\!\!Sp(X,\MM)$.
Namely, consider the function 
$\alpha:TTT(X,\MM)\lra T\!\!{\rm -}\!\!Sp(X,\MM)$, 
defined by $\alpha(\OO)=(X,\BB_\OO,\MM)$ (see Fact~\ref{BO} for the 
notation $\BB_\OO$), and the function
$\beta:T\!\!{\rm -}\!\!Sp(X,\MM)\lra TTT(X,\MM)$, defined by
$\beta\left((X,\BB,\MM)\right)=\OO_\BB$ (see
Corollary~\ref{cord} and Definition~\ref{generates}); then
$\alpha$ and
$\beta$ are bijections and each one is the inverse of the other one.
\end{pro}

\begin{proof}
Let us show that
the function $\alpha$ is well-defined. Let $\OO\in TTT(X,\MM)$. By
Fact~\ref{BO}, the family
$\BB_\OO$ is closed under finite intersections and $X\in\BB_\OO$. By
Proposition~\ref{BO Mclosed}, $\BB_\OO$ is
$\MM$-closed. Hence $(X,\BB_\OO,\MM) \in T\!\!{\rm -}\!\!Sp(X,\MM)$.

We will prove now that
the function $\beta$ is well defined. Let $\BB\subseteq\PP(X)$ be an
$\MM$-closed family, closed
under finite intersections and such that $X\in\BB$. Then, by
Corollary \ref{cord},
$\OO_\BB$ is a topology of
Tychonoff-type on $\MM$.

Proposition~\ref{BOB=B} gives that $\BB_{\OO_\BB}=\BB$, i.e.\
$\alpha\circ\beta={\rm id}_{T\!\!{\rm -}\!\!Sp(X,\MM)}$.

To show that $\beta\circ\alpha={\rm id}_{TTT(X,\MM)}$, let $\OO$ be a
topology of Tychonoff-type on $\MM$.
Then $\beta(\alpha(\OO))=\OO_{\BB_\OO}$. 
Since $\OO$ is a topology of Tychonoff-type, $\left(\BB_\OO\right)^+_\MM$
is a base of $\OO$. 
On the other hand, $\left(\BB_\OO\right)^+_\MM$ is, by definition, a base
of $\OO_{\BB_\OO}$. 
Hence $\OO=\OO_{\BB_\OO}$.
\end{proof}

\begin{defi}\label {HT}
We denote by $\HH\TT$ (Hypertopologies of Tychonoff-type) the category
defined as follows: its objects are all ordered triples
$(X,\MM,\OO)$ where $X$ is a set, $\MM$ is a natural family in $X$ and $\OO$
is
a topology of Tychonoff-type on
$\MM$.
To define the morphisms of $\HH\TT$, let $(X,\MM,\OO)$,
$(X^\prime,\MM^\prime,\OO^\prime)$ be objects of $\HH\TT$
and
$f: X\lra X^\prime$
be a function
between the sets $X$ and $X^\prime$. We will say that
$f$ {\em generates a morphism $f_H$ of} $\HH\TT$ between $(X,\MM,\OO)$ and
$(X^\prime,\MM^\prime,\OO^\prime)$
if
$f(\MM)\subseteq\MM^\prime$ and the induced function on $\MM$,
$f_m:(\MM,\OO)\lra(\MM^\prime,\OO^\prime)$, defined by $f_m(M):= f(M)$
(where
the $M$ on the left-handside is
regarded as an element of $\MM$ and the $M$ on the right-handside is
regarded as a subset of $X$) is
continuous. The morphisms of
$\HH\TT$ are defined to be all $f_H$ generated in this way.
\end{defi}

\begin{rem}\label{rema123}
It is easy to see that not any continuous map between spaces of the
type $(\MM,\OO)$ appears as some $f_H$
(see Definition~\ref{HT} for the notations). Indeed, let $(X,\TT)$
be a discrete space having more than one point; then the constant
function $c:(\CC L(X),\OO_\TT)\lra (\CC L(X),\OO_\TT)$, defined by
$c(F)=X$ for all $F\in\CC L(X)$, is continuous but is not of the type $f_H$
(here $\OO_\TT$ is the classical Tychonoff topology on $\CC L(X)$ (see
Definition~\ref{Tych top})).
\end{rem}

\begin{defi}\label{TH}
We denote by $\TT\HH$ the category defined as follows: 
its objects are all $T$-spaces
$(X,\BB,\MM)$ (see Definition~\ref{T_M-space}).
To define the morphisms of $\TT\HH$,
let $(X,\BB,\MM)$, $(X^\prime,\BB^\prime,\MM^\prime)$ be objects of 
$\TT\HH$ and $f: X\lra X^\prime$ be a function between the sets $X$ and 
$X^\prime$.
We will say that $f$
{\em generates a morphism $f^T$ of} $\TT\HH$ between $(X,\BB,\MM)$ and
$(X^\prime,\BB^\prime,\MM^\prime)$
if $f(\MM)\subseteq\MM^\prime$ and $f^{-1}(\BB^\prime)\subseteq \BB$.
The morphisms of $\TT\HH$ are defined to be all $f^T$ generated in this 
way.
\end{defi}

\begin{lm}\label {lemma1}
Let $X$, $X^\prime$ be sets, $\MM\subseteq\PP(X)$ and
$\MM^\prime\subseteq\PP(X^\prime)$. Let $f:X\lra
X^\prime$ be a function such that $f(\MM)\subseteq\MM^\prime$ and let
$f_m:\MM\lra\MM^\prime$ be defined by
$f_m(M):=f(M)$.
Then, for all $A^\prime\subseteq X^\prime$, we have
$$f_m^{-1}\left((A^\prime)^+_{\MM^\prime}\right)=
\left(f^{-1}(A^\prime)\right)^+_\MM.$$
\end{lm}

\begin{proof}
Let $M\in f_m^{-1}\left((A^\prime)^+_{\MM^\prime}\right)$. Then
$f_m(M)\in(A^\prime)^+_{\MM^\prime}$, i.e.
$f(M)\subseteq A^\prime$. Since $M\subseteq
f^{-1}\left(f(M)\right)\subseteq f^{-1}(A^\prime)$, we obtain
$M\in
\left(f^{-1}(A^\prime)\right)^+_\MM$.

Let $M\in\left(f^{-1}(A^\prime)\right)^+_\MM$. Then $M\subseteq
f^{-1}(A^\prime)$ and hence $f(M)\subseteq
A^\prime$. Therefore $f_m(M)\in (A^\prime)^+_{\MM^\prime}$, i.e.\ 
$M\in f_m^{-1}\left((A^\prime)^+_{\MM^\prime}\right)$.
\end{proof}

\begin{theo}\label {isomorphism}
The categories $\HH\TT$ and $\TT\HH$ are isomorphic.
\end{theo}

\begin{proof}
We define a functor $F:\HH\TT\lra\TT\HH$ as follows: for all $(X,\MM,\OO)\in
\card{\HH\TT}$,
we put
$F\left((X,\MM,\OO)\right):= (X,\BB_\OO,\MM)$
(see Fact \ref{BO} for the notation $\BB_\OO$);
for all morphisms $f_H:(X,\MM,\OO)\lra
(X^\prime,\MM^\prime,\OO^\prime)$, we put
$F(f_H):=f^T$
(see Definitions \ref{HT} and \ref{TH} for the notations $f_H$ and $f^T$).

We define a functor
$G:\TT\HH\lra\HH\TT$ as follows: for all $(X,\BB,\MM)\in \card{\TT\HH}$,
we put
$G\left((X,\BB,\MM)\right):=(X,\MM,\OO_\BB)$
(see Definition \ref{generates} for the notation $\OO_\BB$);
for all morphisms $f^T:(X,\BB,\MM)\lra
(X^\prime,\BB^\prime,\MM^\prime)$, we put
$G(f^T):= f_H$.

Let us check that $F$ and $G$ are well-defined.

Let $\OO$ be a topology of Tychonoff-type on $\MM$.
By Proposition~\ref{char Tych}, the triple $(X,\BB_O,\MM)$ is a $T$-space,
so
that $F(X,\MM,\OO)\in \card{\TT\HH}$.

Let $f_H:(X,\MM,\OO)\lra(X^\prime,\MM^\prime,\OO^\prime)$ be a morphism in
$\HH\TT$.
Then $f(\MM)\subseteq\MM^\prime$ and the induced function on $\MM$,
$f_m:(\MM,\OO)\lra(\MM^\prime,\OO^\prime)$, defined by $f_m(M)= f(M)$, is
continuous.
To check that 
$$f^T:F(X,\MM,\OO)\lra F(X^\prime,\MM^\prime,\OO^\prime)$$ 
(i.e.\ $f^T:(X,\BB_\OO,\MM)\lra (X^\prime,\BB_{\OO^\prime},\MM^\prime$))
is a morphism in $\TT\HH$,
 we need to show that
$f^{-1}(\BB_{\OO^\prime})\subseteq \BB_\OO$.
Let $B^\prime\in\BB_{\OO^\prime}$. Then
$(B^\prime)^+_{\MM^\prime}\in \OO^\prime$.
 By the continuity of $f_m$,
we have that
$f_m^{-1}\left((B^\prime)^+_{\MM\ap}\right)\in\OO$.
By Lemma~\ref{lemma1},
$f_m^{-1}\left((B^\prime)^+_{\MM^\prime}\right)=
\left(f^{-1}(B^\prime)\right)^+_\MM$.
Hence
$\left(f^{-1}(B^\prime)\right)^+_\MM\in \OO$,
i.e.\ $f^{-1}(B^\prime)\in\BB_\OO$.
So, we have proved that $F$ is well-defined. Clearly, $F$ is a functor.

Let now $(X,\BB,\MM)\in\card{\TT\HH}$. 
Then, by Proposition~\ref{char Tych}, the topology $\OO_\BB$ is a
Tychonoff-type topology on $\MM$.
Hence, $G((X,\BB,\MM))\in\card{\HH\TT}$.

Let $f^T:(X,\BB,\MM)\lra(X^\prime,\BB^\prime,\MM^\prime)$ be a morphism in 
$\TT\HH$.
Then $f(\MM)\subseteq\MM^\prime$ and
$f^{-1}(\BB^\prime)\subseteq \BB$.
To check that
$$f_H:(X,\MM,\OO_\BB)\lra (X^\prime,\MM^\prime,\OO_{\BB^\prime})$$ 
is a morphism in $\HH\TT$,
we need to check that the function
$f_m:(\MM,\OO_\BB)\lra(\MM^\prime,\OO_{\BB^\prime})$,
defined by $f_m(M):= f(M)$, is continuous.
By definition (see \ref{generates}),
$(\BB\ap)^+_{\MM\ap}$
is a base of the topology
$\OO_{\BB\ap}$.
Let $B\ap\in\BB\ap$. Then
$(B\ap)^+_{\MM\ap}\in (\BB\ap)^+_{\MM\ap}$.
By assumption,
$f^{-1}(B^\prime)\in \BB$. Hence $\left(f^{-1}(B^\prime)\right)^+_\MM\in
\OO_\BB$. Since,
by Lemma~\ref{lemma1},
$\left(f^{-1}(B^\prime)\right)^+_\MM=
f^{-1}_m\left((B^\prime)^+_{\MM^\prime}\right)$,
we obtain that
$f^{-1}_m\left((B^\prime)^+_{\MM^\prime}\right)\in\OO_\BB$. Therefore, the
function $f_m$ is continuous.
So, $G$ is well-defined. Obviously, $G$ is a functor.

By Proposition~\ref{char Tych}, we have $F\circ G = {\rm id}_{\TT\HH}$ and
$G\circ F = {\rm id}_{\HH\TT}$ on the
objects. The equalities are clearly true for the morphisms. Hence $F$ and
$G$ are isomorphisms.
\end{proof}

We recall that a topological space $(X,\TT)$ is called a
{\it $P_\infty$-space} (see \cite{A, DV}) if
$\TT$ is closed under arbitrary intersections.

\begin{lm}\label{lemma2}
A space $(X,\TT)$ is a $P_\infty$-space if and only it it has a base $\BB$ 
closed under arbitrary intersections.
\end{lm}

\begin{proof}
Assume $\TT$ has a base $\BB$ closed under arbitrary intersections. 
Let $\UU\subseteq\TT$. 
Since $\emptyset$ is an open set, we can assume that 
$\bigcap\UU\not=\emptyset$. 
Let $x\in\bigcap\UU$. 
For any $U\in\UU$, let $B_U\in\BB$ be such that $x\in B_U\subseteq U$. 
Then 
$$x\in\bigcap \{B_U:U\in\UU\}\subseteq\bigcap\UU$$ 
and, by assumption, $\bigcap\{B_U:U\in\UU\}\in\BB$.
Hence, $\bigcap\UU\in\TT$.
\end{proof}

\begin{pro}\label{Pinf}
Let $(X,\MM,\OO)\in\card{\HH\TT}$.
The topological space $(\MM,\OO)$ is
a $P_\infty$-space if and only
if the family $\BB_\OO$ is closed under arbitrary intersections.
If $(\MM,\OO)$ is a $P_\infty$-space then
$(X,\TT_\OO)$ is a $P_\infty$-space.
If $\OO$ is a Tychonoff topology on $\MM$, then $(\MM,\OO)$ is a
$P_\infty$-space if and only if $(X,\TT_\OO)$ is a $P_\infty$-space
(see Fact \ref{BO} for the notations $\BB_\OO$ and $\TT_\OO$).
\end{pro}

\begin{proof}
The first two assertions follow from Lemma \ref{lemma2}, Fact
\ref{intersection}(a) and the definitions of $\BB_\OO$ and $\TT_\OO$.
The last assertion follows now from Corollary~\ref{char Tych top}.
\end{proof}

\begin{cor}\label{Pinf subcat}
Let $\HH\TT_\infty$ be the full subcategory of $\HH\TT$
having as objects all triples
$(X,\MM,\OO)\in \card{\HH\TT}$ for which the
space $(\MM,\OO)$ is a $P_\infty$-space. Let $\TT\HH_\infty$ be the full
subcategory of $\TT\HH$ whose objects are all
$(X,\BB,\MM)\in \card{\TT\HH}$
such that the family $\BB$ is closed under arbitrary intersections.
Then $\HH\TT_\infty$ and $\TT\HH_\infty$ are isomorphic.
\end{cor}

\begin{proof}
It follows from (the proof of) Theorem~\ref{isomorphism} and
Proposition~\ref{Pinf}.
\end{proof}

\begin{exa}\label{exaPinf}
We will show that there exists $(X,\MM,\OO)\in\card{\HH\TT}$
such that $(X,\TT_\OO)$ is a $P_\infty$-space but $(\MM,\OO)$
is not a $P_\infty$-space.

Let $X=\omega$,
$$
\BB=\{\{n\}:n\in\omega\}\cup
\{A\subseteq \omega : |\omega\setminus A|<\aleph_0\}\cup \{\ems\},
$$
and $\PP(\omega)\spe\MM\spe\FF in_2 (\omega)\stm \{\ems\}$.
 The
family $\BB$ is a base of the discrete topology $\TT$ on $\omega$, it is
closed
under finite
intersections (but not under infinite intersections) and $X\in\BB$. Let us
show that $\BB$
is $\MM$-closed.

Let $\{ B_\d\}_{\d\in\Delta}$ be a subfamily of $\BB\stm\{\ems\}$ which is
an $\MM$-cover of a subset $B$ of $X$.
Without loss of generality, we can assume that there exist
at least two indices $\delta_1$ and $\delta_2$ such that
$B_{\d_1}\not=B_{\d_2}$. Then there is at least one
$\delta\in\Delta$ such that
$|\omega\setminus B_\delta|<\aleph_0$; otherwise we would have
$B_{\d_1}=\{n_{\d_1}\}$,
$B_{\d_2}=\{n_{\d_2}\}$ and the set
$F=\{n_{\d_1},n_{d_2}\}$, which belongs to $\MM$,
would be contained in $B$ without being contained
in $B_\delta$ for
any
$\delta\in\Delta$. Therefore $|\omega\setminus B|<\aleph_0$ and, hence,
$B\in\MM$.

Put $\OO:=\OO_\BB$. Then, by Proposition~\ref{BOB=B}, $\BB=\BB_\OO$
and hence $\TT=\TT_\OO$. Since $\BB$ is not closed under arbitrary
intersections, we obtain, by Proposition~\ref{Pinf}, that $(\MM,\OO)$
is not a $P_\infty$-space.
Clearly,
the space $(X,\TT)$ is a $P_\infty$-space because
it is discrete. Observe that $\OO$ is not a Tychonoff topology since
$\BB_\OO\not=\TT_\OO$.
\end{exa}

\begin{exa}\label{example}
Two more examples of Tychonoff-type, non Tychonoff topologies
on some families $\MM\sbe\PP(X)$.

Let $X$ be a set with more than two elements. Let $\MM=\FF
in(X)\setminus\{\emptyset\}$ and let
$\OO=\{\{\{x\}:x\in A\} : A\subseteq X\}\cup\{\MM\}\cup\{\emptyset\}$. Then
$\OO$ is a topology on $\MM$.

$\OO$ is a topology of Tychonoff-type since 
$$
\OO\cap\PP(X)^+_\MM=\{ \{x\} : x\in X\}\cup\{X^+_\MM\}\cup\{\emptyset\}$$ 
is a base for $\OO$.
Clearly $(\MM,\OO)$ is a $P_\infty$-space.

We have $\BB_\OO=\{\{x\}: x\in X\}\cup\{X\}\cup\{\emptyset\}$, and
therefore $\TT_\OO$ is the discrete
topology. By Proposition~\ref{char Tych top}, since $\BB_\OO\not=\TT_\OO$,
$\OO$
is not a Tychonoff topology on
$\MM$.

Observe that $(\MM,\OO)$ is not a T$_0$-space. In fact, the only
neighbourhood
of any $F\in\MM$ such that
$|F|\geq 2$ is $\MM$, and we are assuming $|X|>2$.

We consider now the natural family $\MM^\prime=\PP(X)\setminus\{\emptyset\}$
and we define $\OO$ as above.
$\OO$ is a Tychonoff-type topology on $\MM^\prime$ but it is not a Tychonoff
topology.
Again $\TT_\OO$ is the discrete topology on $X$. Hence $\MM^\prime = \CC
L((X,\TT_\OO))$.
We observe as before that $(\MM\ap,\OO)$ is not a T$_0$-space.
Note also that the family $\BB_\OO$ is $\MM^{\prime\prime}$-closed for every
natural family
$\MM^{\prime\prime}$ in $X$ which contains all two-points subsets of $X$ and
$\emptyset\not\in\MM^{\prime\prime}$.
\end{exa}

We will briefly discuss now some topological
properties of the hyperspaces $(\MM,\OO)$ with Tychonoff-type topologies
$\OO$.

\begin{fact}\label{sep axioms}
Let $X$ be a set, $\MM\subseteq\PP(X)$ and $\OO$ be a topology of
Tychonoff-type on $\MM$ generated by a subfamily $\BB$ of $\PP(X)$.
Then:
\begin{itemize}
\item[(a)]
the topological space $(\MM,\OO)$ is
a $T_0$-space (resp., $T_1$-space) if and only if for any
$F,G\in\MM$ with $F\not= G$, there exists a $B\in\BB$ such that either
$F\subseteq B$ and
$G\not\subseteq B$, or $G\subseteq B$ and $F\not\subseteq B$ (resp.,
$F\subseteq B$ and
$G\not\subseteq B$).
\item[(b)]
if for any $x\in X$ and for any $F\in\MM$ with $x\not\in F$,
there exists
a $B\in\BB$ such that $F\subseteq B$ and $x\not\in B$, then $(\MM,\OO)$
is a $T_0$-space.
\end{itemize}
\end{fact}

\begin{rem}\label{Choban}
Let us note that Fact~\ref{sep axioms}(b) implies the following assertion,
which was mentioned in \cite{C}, section 2 (after Lemma 3)
(the requirement that $X\in\Omega$ has
to be added there):
{\it if $(X,\TT)$ is a regular $T_1$-space, $\MM$ is a family consisting of
closed subsets of $(X,\TT)$ and $\BB$ is a base of $(X,\TT)$ such that
$\BB=\BB^\cap$ and
$U\in\BB$ implies that $X\stm \ovl{U}\in\BB$, then $(\MM,\OO_\BB)$ is
a $T_0$-space.}
\end{rem}

\begin{fact}\label{T_0}
Let $(X,\MM,\OO)\in\card{\HH\TT}$. Then the correspondence
$(X,\TT_\OO)\lra (\MM,\OO)$, $x\mapsto \{x\}$, is a homeomorphic embedding.
Hence, we have, in particular, that:
\begin{itemize}
\item[(a)] 
$w(X,\TT_\OO)\le w(\MM,\OO)$;
\item[(b)]
if $(\MM,\OO)$ is a $T_0$-space then $(X,\TT_\OO)$ is a $T_0$-space.
\end{itemize}
\end{fact}

\begin{fact}\label{FsubsetG}
\mbox{}
\begin{itemize}
\item[(a)]
Let $X$ be a set, $\MM\subseteq\PP(X)$ be a family such that there exist
$F, G\in\MM$ with $F\subset G$ and $F\not= G$, and let $\OO$ be a topology
of Tychonoff-type on $\MM$. Then $(\MM,\OO)$ is not a $T_1$-space.
\item[(b)] 
If $(X,\MM,\OO)\in\card{\HH\TT}$ then $(\MM,\OO)$ is a $T_1$-space
if and only if $(X,\TT_\OO)$ is a $T_1$-space and $\MM=\{\{x\}: x\in X\}$.
\end{itemize}
\end{fact}

\begin{fact}\label{compactness}
Let $(X,\MM,\OO)\in\card{\HH\TT}$.
Then $(\MM,\OO)$ is a compact space if and only if any $\MM$-cover of $X$,
consisting of elements of $\BB_\OO$, has a finite $\MM$-subcover.
\end{fact}

\begin{proof}
It follows from Proposition~\ref{union}.
\end{proof}

\begin{exas}
There are many examples of ``very nice'' spaces $X$ with 
non-$T_0$-hyperspa\-ces $(\MM,\OO_\BB)$ (see Examples~\ref{example} and
\ref{bigexample}).
As an example of a non-$T_0$-space $(X,\TT)$ with a $T_0$-hyperspace 
$(\MM,\OO_\TT)$, consider the two-points space $X=\{0,1\}$, with
$\TT=\MM=\{\emptyset,X\}$.

There exist non-compact spaces $X$ such that $(\CC L (X),\OO_\BB)$
is a compact non-$T_0$-space (e.g., the space $(\CC L (\RRRR),\OO_\BB)$,
described in Example~\ref{bigexample}).

To get an example of a non-compact space $X$ and a natural family $\MM$
such that $(\MM,\OO_\BB)$ is a compact $T_0$-space, consider $X:=\RRRR$
with its natural topology, $\MM:=\FF in_2(\RRRR) \cup\{\RRRR\}$ and take
$\BB$ as in Example~\ref{bigexample}.

As an example of a compact space $(X,\TT)$ with a non-compact hyperspace
$(\MM,\OO_\TT)$,
consider the unit interval $X=[0,1]$ with its natural topology and put
$\MM=\{\{x\}:x\in(0,1]\}$.
\end{exas}

The next three propositions are generalizations of, respectively,
Propositions 1, 2 and 3 of \cite{F}, and have proofs similar to those
given in \cite{F}. 
(Let us note that in Proposition 2 of \cite{F} the requirement
``$\ems\not\in\CC$'' has to be added.)

\begin{pro}\label{F1}
Let $(X,\MM,\OO)\in\card{\HH\TT}$, 
$w(\MM,\OO)=\aleph_0$, 
$(X,\TT_\OO)$ be a $T_1$-space,
$\BB_\OO$ be closed under countable unions and
$\MM$ contain all infinite countable closed subsets of $(X,\TT_\OO)$.
Then $(X,\TT_\OO)$ is a compact space.
\end{pro}

\begin{pro}\label{F2}
Let $(X,\MM,\OO)\in\card{\HH\TT}$ and $\ems\not\in\MM$. Then
$d(\MM,\OO)=d(X,\TT_\OO)$.
\end{pro}

\begin{pro}\label{F3}
Let $(X,\MM,\OO)\in\card{\HH\TT}$ and $\ems\not\in\MM$. Then
$(\MM,\OO)$ has isolated points if and only if $(X,\TT_\OO)$
has isolated points.
\end{pro}

\section{On $\OO$-commutative spaces}

\begin{nist}\label{2-emb}
Let $(X,\TT)$ be a topological space and $A\subseteq X$. Recall that
$A$ is
said to be
{\em 2-combinatorially embedded in} $X$ (see \cite{CN}) if
the closures in $X$ of any two disjoint closed in $A$
subsets of $A$ are disjoint.
\end{nist}

\begin{defi}\label{2Bcomb}
Let $(X,\TT)$ be a topological space, $A\sbe X$
and $\BB\subseteq \PP(X)$.
We will
say that $A$ is {\em $2_\BB$-combinatorially embedded in $X$} if for any
$F\in \CC L(A)$ and for any $U\in\BB$
with $F\subseteq U$, there exists a $V\in\BB$ such that $\overline{F}^X
\subseteq V$ and $V\cap A\subseteq U$.
\end{defi}

\begin{pro}\label{2T}
Let $(X,\TT)$ be a topological space and $A\subseteq X$. Then $A$ is
2-com\-bi\-na\-to\-rially embedded in
$X$ if and only if $A$ is $2_\TT$-combinatorially embedded in $X$.
\end{pro}

\begin{proof}
($\Ra$)
Let $H\in \CC L(A)$, $V\in\TT$ and $H\sbe V$.
We put $U=V\cap A$ and
$F=A\setminus U$. Then $F$ and $H$ are two disjoint closed subsets of $A$.
Hence, by assumption, they have disjoint closures in $X$, i.e.\ 
$\overline{F}^X\cap\overline{H}^X=\emptyset$.
Let
$W=X\setminus\overline{F}^X$.
Then $W$ is open in $X$, $W\cap A=U=V\cap A$ and $\overline{H}^X\subseteq
W$.

($\La$)
Let $F$ and $G$ be two disjoint closed subsets of $A$. Put
$V=X\setminus\overline{G}^X$. Then $V$ is
open in $X$ and $F\subseteq V$. Hence, by assumption, there exists an open
set $W$ such that $\overline{F}^X\subseteq W$
and $W\cap A\subseteq V\cap A$. Let $U=A\setminus G$.
Then $Ex_{A,X}U=V$.
Hence $V\cap A=U$ and $W\subseteq V$.
We conclude that $\overline{F}^X\subseteq
W\subseteq
V=X\setminus\overline{G}^X$, i.e.
$\overline{F}^X\cap\overline{G}^X=\emptyset$.
\end{proof}

\begin{rem}\label{rem111}
In Example~\ref{exa111} below
we will show that there exist spaces $(X,\TT)$, subspaces
$A$ of $X$ and bases $\BB$ of $\TT$ such that $A$ is
$2_\BB$-combinatorially embedded in $X$
but $A$ is not
2-combinatorially embedded in $X$.
\end{rem}

\begin{pro}\label{iAX}
Let $(X,\TT)$ be a T$_1$-space, $\OO$ be a topology of Tychonoff-type on
$\CC L(X)$ and
$A\subseteq X$. Put $\BB_A=\{U\cap A : U\in\BB_\OO\}$
(see Fact \ref{BO} for the notation $\BB_\OO$).
The family $\BB_A$ generates a topology of Tychonoff-type $\OO_A$ on $\CC
L(A)$.
The function $i_{A,X}:(\CC L(A),\OO_A)\lra (\CC L(X),\OO)$, defined by
$i_{A,X}(F):=\overline{F}^X$, is
inversely continuous (i.e.\ it is injective and its inverse, defined on
$i_{A,X}\left(\CC L(A)\right)$, is a
continuous function) if and only if the set
$A$ is $2_{\BB_\OO}$-combinatorially embedded in $X$.
\end{pro}

\begin{proof}
The family $\BB:=\BB_\OO$ is closed under
finite intersections and $X\in\BB$ (see Fact \ref{BO}).
Hence
the family $\BB_A$ is closed under finite intersections and $A\in\BB_A$.
Therefore, by Corollary \ref{cord},
$\BB_A$ generates a topology of Tychonoff-type $\OO_A$ on $\CC L(A)$.

The function $i_{A,X}:(\CC L(A),\OO_A)\lra (\CC L(X),\OO)$ is clearly
injective.
Denote by $g$ its inverse defined on
$i_{A,X}(\CC L(A))$, i.e.
$g:i_{A,X}(\CC L(A))\lra \CC L(A)$.

($\Ra$)
Let $H\in\CC L(A)$, $U\in\BB$ and $H\subseteq U$. Then $H\in
\left(U\cap A\right)^+_{\CC L(A)}\in \OO_A$. Since
$g(\ovl{H}^X)=H$, the continuity of $g$ implies that
there exists a $V\in \BB$ such that
$$\ovl{H}^X\sbe V
\mbox{ and }
g(V^+_{\CC L(X)}\cap i_{A,X}(\CC L(A)))\sbe (U\cap A)^+_{\CC L(A)}.$$ 
Then
$V\cap A\subseteq U\cap A$.
Indeed, let $x\in V\cap A$.
Since $X$ is a $T_1$-space, we obtain that
$$
\{x\}\in V^+_{\CC L(X)}\cap i_{A,X}(\CC L(A))
\mbox{ and }
g(\{x\})=\{x\}.$$
Hence $x\in U\cap A$. So,
$A$ is $2_\BB$-combinatorially embedded in $X$.

($\La$)
Let $F\in i_{A,X}(\CC L(A))$
and $g(F)=H$. Then $F=\ovl{H}^X$ and
$H\in \CC L(A)$. Let $U\in\BB_A$ be such that $H\sbe U$. Then there
exists a $V\in \BB$ with $V\cap A=U$. Hence $H\sbe V$. Since
$A$ is $2_\BB$-combinatorially embedded in $X$,
there exists a $W\in\BB$ such that $F=\ovl{H}^X\sbe W$ and
$W\cap A\sbe V\cap A=U$.
Then $F\in W^+_{\CC L(X)}\in\OO$. 
We will show that
$$g(W^+_{\CC L(X)}\cap i_{A,X}(\CC L(A)))\sbe U^+_{\CC L(A)}.$$
Indeed,
let $K\in \CC L(A)$,
$G=\ovl{K}^X$
and $G\sbe W$. Then $g(G)=K$
and 
$$K=G\cap A\sbe W\cap A\sbe V\cap A=U,$$
i.e.\ $K\in U^+_{\CC L(A)}$, as we have to show. Hence, $g$ is a
continuous function.
\end{proof}

\begin{cor}[\cite{D2}, Theorem 2.1]\label{2.1D2}
If in Proposition \ref{iAX} we take $\OO$ to be the
Tychonoff topology on
$\CC L(X)$
generated by $(X,\TT)$
then
the function $i_{A,X}$ is inversely continuous if and only if
$A$ is 2-combinatorially embedded in $X$.
\end{cor}

\begin{proof}
It follows from Propositions~\ref{iAX}, \ref{uniqueness} and \ref{2T}.
\end{proof}

\begin{cor}\label{sequential}
Let $(X,\TT)$ be a T$_2$-space, $A\subseteq X$ and $\OO$ be a topology of
Tychonoff-type on $\CC L(X)$
generated by a subfamily of $\TT$.
Let $i_{A,X}$ be inversely continuous
(see Proposition~\ref{iAX} for the notation $i_{A,X}$).
Assume that the following condition is satisfied:
\begin{itemize}
\item[(*)] For any $U\in\TT$ and for
all countable $F\in\CC L(A)$ such that $\card{A\setminus
F}\geq\aleph_0$ and
$F\subseteq U$, there exists a $V\in\BB_\OO$
with $F\subseteq V\subseteq U$.
\end{itemize}
Then the set $A$ is sequentially closed.
\end{cor}

\begin{proof}Put $\BB :=\BB_\OO$.
Then, by Proposition~\ref{char BOB}, $\BB\sbe \TT$.
Assume that the set $A$ is not sequentially closed.
Then there exists a sequence
$(x_n)_{n\in\omega}$ in $A$ and a point $x\in X\setminus A$ such that
$\lim_{n\lra\infty}x_n=x$. Without loss
of generality we can assume $x_n\not=x_m$ for all $n\not=m$.

Let us consider the sets $F=\{x_{2n} : n\in\omega\}$ and
$G=\{x_{2n-1} : n\in\omega\}$. 
Put $U=X\setminus \overline{G}^X$. 
Then $F$ is a countable closed subset of $A$,
$\card{A\stm F}\ge\aleph_0$, $F\subseteq U$ and $U\in \TT$.
By (*), there exists a $V\in\BB$ such that $F\subseteq V\subseteq U$. 
Since we are assuming that the function $i_{A,X}$ is inversely continuous, 
we obtain, by Proposition~\ref{iAX}, that the set $A$ is 
$2_\BB$-combinatorially embedded in $X$. 
Hence there exists a $W\in\BB$ such that $\overline{F}^X\subseteq W$ and 
$W\cap A\subseteq V\cap A$.
Then $x\in W$, because $x\in\ovl{F}^X$.
Since $W\in\TT$ and $x$ is a limit point of $G$, we have 
$G\cap W\not=\emptyset$. 
However this is a contradiction because 
$$W\cap A\subseteq V\cap A\sbe U=X\setminus\overline{G}^X,$$
and hence $G\cap W=\ems$.
Therefore, $A$ is sequentially closed.
\end{proof}

\begin{rem}\label{rem222}
In Example~\ref{exa111} below we will show that condition (*) of 
Corollary~\ref{sequential} is essential, i.e., if we omit it, then the set 
$A$ could fail to be sequentially closed.
\end{rem}

\begin{cor}[\cite{D2}, Corollary 2.3]\label{2.3D2}
Let $(X,\TT)$ be a T$_2$-space, $A\subseteq X$, $\OO$ be the
Tychonoff topology on $\CC L(X)$
generated by $(X,\TT)$
and $i_{A,X}$ be inversely continuous
(see Proposition~\ref{iAX} for the notation $i_{A,X}$).
Then the set $A$ is sequentially closed.
\end{cor}

\begin{proof}
We have, by Proposition~\ref{uniqueness}, that $\BB_\OO=\TT$.
Then condition (*) of
Corollary \ref{sequential}
is trivially satisfied.
Hence, by
Corollary \ref{sequential},
$A$ is sequentially closed.
\end{proof}

\begin{cor}\label{eqv}
Let $(X,\TT)$ be a sequential T$_2$-space, $A\subseteq X$ and $\OO$ be
a topology of
Tychonoff-type on $\CC L(X)$
generated by a subfamily of $\TT$.
Assume that condition (*) of Corollary~\ref{sequential} is satisfied.
Then the following conditions are
equivalent
(see Proposition~\ref{iAX} for the notation $i_{A,X}$):
\begin{itemize}
\item[(a)] $i_{A,X}$ is a homeomorphic embedding;
\item[(b)] $i_{A,X}$ is inversely continuous;
\item[(c)] $A$ is closed in $X$.
\end{itemize}
\end{cor}

\begin{proof}
It is clear that (a) implies (b).
The implication (c)$\Ra$(a) is true for any $X$, because
if $A$ is a
closed subset of $X$ then $i_{A,X}$ is the inclusion map.
Let us show that (b) implies (c). By Corollary~\ref{sequential}, $A$ is
sequentially closed. Since $X$ is a
sequential space, we obtain that the set $A$ is closed.
\end{proof}

\begin{cor}[\cite{D2}, Corollary 2.4]\label{2.4D2}
Let $(X,\TT)$ be a sequential T$_2$-space, $A\subseteq X$ and $\OO$ be the
Tychonoff topology on $\CC L(X)$ generated by $(X,\TT)$.
Then the following conditions are equivalent
(see Proposition~\ref{iAX} for the notation $i_{A,X}$):
\begin{itemize}
\item[(a)] $i_{A,X}$ is a homeomorphic embedding;
\item[(b)] $i_{A,X}$ is inversely continuous;
\item[(c)] $A$ is closed in $X$.
\end{itemize}
\end{cor}

\begin{defi}\label{O-commutative}
Let $(X,\TT)$ be a topological space and let $\OO$ be a topology of 
Tychonoff-type on $\CC L(X)$.
The space $(X,\TT)$ is called {\em $\OO$-commutative } if for any 
$A\subseteq X$ the function $i_{A,X}$, defined in Proposition~\ref{iAX},
is a homeomorphic embedding.

When $\OO$ is the Tychonoff topology on $\CC L(X)$ generated by $(X,\TT)$,
the notion of ``$\OO$-commutative space'' coincides with the notion of
``{\em commutative space}'', introduced in \cite{D1, D2}.
\end{defi}

\begin{cor}\label{comm}
Let $(X,\TT)$ be a sequential T$_2$-space, $\OO$ be
a topology of
Tychonoff-type on $\CC L(X)$
generated by a subfamily of $\TT$ and
let condition (*) of Corollary~\ref{sequential} be satisfied
for every subspace $A$ of $X$.
Then $X$ is $\OO$-commutative
if and only if
$X$ is discrete.
\end{cor}

\begin{proof}
It follows from Corollary~\ref{eqv}.
\end{proof}

\begin{cor}[\cite{D2}, Corollary 2.5]
\label{2.5D2}
If $X$ is a sequential T$_2$-space then $X$ is commutative if and only if
$X$ is discrete.
\end{cor}

\begin{exa}\label{ex F}
Let us show that there exist spaces $X$ and
topologies $\OO$ of Tycho\-noff-type on $\CC L(X)$ that are
not Tychonoff topologies and that satisfy all hypothesis of Corollary
\ref{comm}.

Let $X=D(\aleph_1)$ be the discrete space of cardinality $\aleph_1$.
Let 
$$\BB=\{A\subset X : |A|\leq\aleph_0\}\cup\{X\}$$ 
and $\MM=\CC L(X)$. 
Then $\MM$ is a natural family on $X$, $\BB$ is $\MM$-closed,
$\BB^\cap=\BB$, $X\in\BB$ and
$\BB$ is a base for the discrete topology on $X$.
Let $\OO_\BB$ be the topology on $\MM$ generated
by $\BB$. Then $\OO_\BB$ is a
topology of Tychonoff-type on $\MM$,
however it is not a Tychonoff topology. In fact, by Proposition~\ref{BOB=B},
$\BB=\BB_{\OO_\BB}$. Hence
$\BB_{\OO_\BB}\not=\TT_{\OO_\BB}$ and,
by Corollary~\ref{unique}, $\OO$ cannot be a
Tychonoff topology.
Obviously, $\BB=\BB_{\OO_\BB}$ satisfies condition (*) of
Corollary \ref{sequential} for any subspace $A$ of $X$.
\end{exa}

\begin{defi}\label{O-HS-space}
Let $(X,\TT)$ be a topological space and $\OO$ be a topology of
Tycho\-noff-type on $\CC L(X)$.
The space $(X,\TT)$ is called {\em $\OO$-HS-space} if, for any
$A\subseteq X$, the function $i_{A,X}$,
defined in Proposition~\ref{iAX},
is continuous.

When
$\OO$ is the Tychonoff topology on $\CC L(X)$
generated by $(X,\TT)$,
the notion of ``$\OO$-HS-space'' coincides with the notion of
``{\em HS-space}'', introduced in \cite{BDN1, BDN2}.
\end{defi}

\begin{cor}\label{ocom}
Let $(X,\TT)$ be a T$_1$-space and $\OO$ be a topology of Tychonoff-type on
$\CC L(X)$.
Then $X$ is an $\OO$-commutative space if and only if $X$ is an
$\OO$-HS-space
and every subset $A$ of $X$ is $2_{\BB_\OO}$-combinatorially embedded in
$X$.
\end{cor}

\begin{proof}
It follows from Proposition~\ref{iAX}.
\end{proof}

\begin{cor}[\cite{D2},Corollary 2.2]\label{ocom1}
A $T_1$-space $X$ is commutative
if and only if $X$ is an HS-space
and every subspace of $X$ is 2-combinatorially embedded in $X$.
\end{cor}

\begin{proof}
It follows from Corollary~\ref{ocom} and Proposition~\ref{2T}.
\end{proof}

\begin{exa}\label{bigexample}
We will describe two Tychonoff-type, non Tychonoff topologies on two
different subfamilies of
$\PP(\RRRR)$ generated by the family $\BB$
of all open intervals of $\RRRR$. One of the resulting spaces
will be $T_0$ and
the other one will not.

Let $\TT$ be the natural topology on $X :=\RRRR$.
Then the family $\BB$ of all
open intervals in
$X$ is a base for $\TT$, it is closed under finite intersections
and $X\in\BB$.
Put $\MM :=\CC L(X,\TT)$ and
$\MM^\prime :=\FF in_2(X)$.
They are natural families.
The family $\BB$ is both $\MM^\prime$-closed and $\MM$-closed.
Indeed,
let $U\subseteq X$ be $\MM^\prime$-covered by a subfamily $\BB_U$ of $\BB$.
Then $U\in\TT$ and for every $x,y\in U$ there exists an open interval
$(\a,\b)\in\BB_U$ containing the points $x$ and $y$.
Hence $U$ is a connected open set in
$\RRRR$, i.e.\ $U\in\BB$.
Therefore, $\BB$ is an $\MM\ap$-closed family.
Since
$\MM^\prime\subset
\MM$, we obtain, by Proposition~\ref{M and M^prime}, that
$\BB$ is an $\MM$-closed family as well.

By Corollary~\ref{cord}, $\BB$ generates Tychonoff-type topologies
$\OO_\BB$ on $\MM$ and
$\OO^\prime_\BB$ on $\MM^\prime$.
As it follows from
Proposition~\ref{BOB=B},
$\BB_{\OO_\BB}=
\BB_{\OO^\prime_\BB}=
\BB\not=\TT$. Hence, by
Corollary~\ref{unique},
$\OO_\BB$
and $\OO^\prime_\BB$
are not Tychonoff topologies on $\MM$, respectively $\MM\ap$.

It is easy to see that $(\MM^\prime,\OO^\prime_\BB)$ is a
T$_0$-space. Indeed, let $\{x,y\}$ and $\{u,v\}$ be two distinct elements
in $\MM^\prime$. 
We can assume $x<x+\varepsilon <u\leq v$ for some 
$\varepsilon>0$. 
Consider the interval $B=(x+\varepsilon,+\infty)$.
Then $\{u,v\}\in B^+_{\MM^\prime}$ but 
$\{x,y\} \not\in B^+_{\MM^\prime}$.

Let's prove that $(\MM,\OO_\BB)$ is not a T$_0$-space.
Put $F=\{2k:k\in\Z\}$ and $G=\{2k+1:k\in\Z\}$. Then
$F,G\in\MM$ and $F\not=G$ but the only neighbourhood of
both $F$ and $G$ in
$\MM$ is $X^+_{\MM}=\MM$.
\end{exa}

\begin{exa}\label{bigexample1}
In the notations of Example~\ref{bigexample}, we will show that
$(\RRRR,\TT)$ {\em is an $\OO_\BB$-HS-space}.

We are working now with the space $(\MM,\OO_\BB)$ from
Example~\ref{bigexample}.
We will write simply $\OO$ instead of $\OO_\BB$.

Let $A\subseteq X$. We have to show that the function 
$$i_{A,X}:(\CC L(A),\OO_A)\lra (\CC L(X),\OO),$$
where the topology $\OO_A$ on $\CC L(A)$ is generated by the family
$$\BB_A=\{A\cap U: U\in\BB\},$$
is continuous
(see Proposition~\ref{iAX} for the notation $i_{A,X}$).
Let
$B\in\BB$.
We will show that
$i_{A,X}^{-1}(B^+_\MM)$ is an open set.
Take an $F\in i_{A,X}^{-1}(B^+_\MM)$. Then $F\in\CC L(A)$ and
$\overline{F}^X\subseteq B$.
There exists an
$E\in\BB$ such that
$$\overline{F}^X\subseteq E\subseteq \overline{E}^X\subseteq B$$
(this is clear if $F$ is bounded, since in this case $\overline{F}^X$ is 
compact; if $F$ is unbounded below, but is bounded above, then 
$B=(-\infty,\beta)$, for some $\b\in\RRRR$, and we can pick 
$E=(-\infty,\gamma)$ with $\sup F<\gamma<\beta$; similarly if $F$ is 
unbounded above but not below; if $F$ is unbounded both above and below
then we have $B=\RRRR$ and we put $E := B$).
Then
$$F\in \left( E\cap A\right)^+_{\CC L(A)}\subseteq 
i_{A,X}^{-1}(B^+_\MM).$$ 
Indeed, let $G\in \left( E\cap A\right)^+_{\CC L(A)}$. 
Then 
$$\overline{G}^X\subseteq \overline{E}^X\subseteq B,$$
i.e.\ $i_{A,X}(G)\in B^+_\MM$.
\end{exa}

\begin{rem}\label{rem112}
Let us note that
a similar proof shows that
{\em every subspace $Y$ of $(\RRRR,\TT)$
is an
$\OO_{\BB_Y}$-HS-space}
(see Examples~\ref{bigexample} and
\ref{bigexample1} for the notations).

More generally, let $Y$ be a topological space and $\DD$ be a base of $Y$.
We will say that $Y$ is {\it $\DD$-normal} if for every $F\in\CC L(Y)$
and for every $D\in\DD$ such that $F\sbe D$ there exists an $E\in\DD$ with
$F\sbe E\sbe \ovl{E}^Y\sbe D$.
Now, arguing as in Example~\ref{bigexample1},
we can prove that
{\em if $Y$ is a $\DD$-normal space, $\DD=\DD^\cap$ and $Y\in\DD$,
then $Y$ is an $\OO$-HS-space, where $\OO$ is the Tychonoff-type
topology on $\CC L(Y)$ generated by $\DD$.}
This generalizes the result of M. Sekanina \cite{Se} that any normal
space is a HS-space.
\end{rem}

\begin{exa}\label{exa111}
In the notations of Examples~\ref{bigexample} and \ref{bigexample1}, we 
will show that {\em the function $i_{A,\RRRR}$ is a homeomorphic embedding
for any open interval $A$.}

We will argue for $A=(0,1)$;
the proof for any other open interval is similar.
We know, by Example~\ref{bigexample1}, that the function $i_{A,X}$ is 
continuous. 
Therefore we only need to prove that $i_{A,X}$ is inversely continuous. 
By Proposition~\ref{iAX}, it is enough to show that the set $A$ is
$2_\BB$-combinatorially embedded in $X$.
So, let $H$ be a closed subset of $(0,1)$ and let $B=(\a,\b)\in\BB$ be 
such that $H\subseteq B$.
We have to find a $D\in\BB$ such that $\overline{H}^X\subset D$ and 
$D\cap (0,1)\subseteq B$.
Clearly, $\overline{H}^X\subseteq [0,1]$. 
If $\overline{H}^X\subset (0,1)$, we can take $D=B$.
If $0\in\overline{H}^X$ but $1\not\in\overline{H}^X$ then $\a\le 0$ and we
can put $D=(-1,\b)$.
If $1\in\overline{H}^X$ but $0\not\in\overline{H}^X$ then $\b\ge 1$ and we
can put $D=(\a,2)$.
If $0,1\in\overline{H}^X$ then $\a\le 0$ and $\b\ge 1$ and we put 
$D=(-1,2)$.
Therefore, {\em $A$ is $2_\BB$-combinatorially embedded in $(\RRRR,\TT)$.}

Note that $A$ {\em is not 2-combinatorially embedded in $(\RRRR,\TT)$.}

Observe that the triple $((\RRRR,\TT), A,\OO)$ satisfies all hypothesis of
Corollary~\ref{sequential} except for condition (*), but $A$
is not sequentially closed. 
\end{exa}

\begin{exa}\label{prop111}
Let $Y\sbe \RRRR$. We will say, as usual, that a point $x\in Y$ is {\em
isolated from the right (left) (in $Y$)} if there exists an $\ep >0$ such
that if we put $U=(x,x+\varepsilon)$ ($U=(x-\varepsilon,x)$) then $U\cap
Y=\emptyset$. Now, in the notations of Examples~\ref{bigexample} and
\ref{bigexample1}, we have: {\em a subspace $Y$ of $(\RRRR,\TT)$ is
$\OO_{\BB_Y}$-commutative if and only if every point of $Y$ is either
isolated from the right or from the left.}

We first show that a space $Y$ that has a point $y_0$ which is
non-isolated both from the left and from the right cannot be
$\OO_{\BB_Y}$-commutative. Indeed, put $A=Y\setminus\{y_0\}$. We will
prove that $A$ is not $2_{\BB_Y}$-combinatorially embedded in $Y$. By
Proposition~\ref{iAX}, this will imply that the function $i_{A,Y}$ is not
inversely continuous and hence the space $Y$ will be not
$\OO_{\BB_Y}$-commutative. Let $H=\{y_n:n\in\omega\}$ be a decreasing
sequence in $Y$ converging to $y_0$. Then $H$ is a closed subset of $A$
and $H\subset (y_0,+\infty)\cap Y$. Suppose that there exists $B\in\BB$
such that $cl_Y H\subseteq B$ and $B\cap A\sbe (y_0,+\infty)$. Since
$y_0\in cl_Y H\sbe B$ and $y_0$ is not isolated from the left, we have
that $(B\cap A)\setminus (y_0,+\infty)\not =\ems$, which is a
contradiction. Hence, $A$ is not $2_{\BB_Y}$-combinatorially embedded in
$Y$.

Now we will show that a space $Y$ having only points which are isolated
either from the left or from the right is $\OO_{\BB_Y}$-commutative. Let
$A\subset Y$. We know, by Remark~\ref{rem112}, that the function $i_{A,Y}$
is continuous. Hence it is enough to show that it is inversely continuous,
i.e., according to Proposition~\ref{iAX}, that $A$ is
$2_{\BB_Y}$-combinatorially embedded in $Y$. So, let $H\in\CC L(A)$ and
let $H\subseteq B\cap Y$ for some $B=(\alpha,\beta)$. We have to find a
$D\in\BB$ such that $\overline{H}^Y\subset D\cap Y$ and $D\cap A\subseteq
B$. We have $ cl_Y H\subseteq \overline{B}^X=[\alpha,\beta]$. If $cl_Y
H\subseteq B$, we can take $D=B$ and we are done. If $\alpha\in cl_Y H$
and $\b\not\in cl_Y H$ then $\alpha\in Y$ and $\alpha$ is not isolated
from the right, being a limit point of $H$. Hence, by the assumption,
$\alpha$ is isolated from the left. Thus there exists a $\gamma<\alpha$
such that $(\gamma,\alpha)\cap Y=\emptyset$. Then $D=(\gamma,\b)$ is the
required interval. The other two possible cases are treated analogously.
\end{exa}


\begin{thebibliography}{10}

\bibitem{A}
P.~S. Alexandroff, \emph{{D}iskrete {R}\"aume}, Rec. Math. [Mat. Sbornik] N.S.
  \textbf{2} (1937), 501--518.

\bibitem{BDN1}
S.~Barov, G.~Dimov, and S.~Nedev, \emph{On a theorem of {H}.-{J}.\ {S}chmidt},
  C. R. Acad. Bulgare Sci. \textbf{46} (1993), no.~3, 9--11. \MR{1 261 974}

\bibitem{BDN2}
\bysame, \emph{On a question of {M}. {P}aoli and {E}. {R}ipoli}, Boll. Un. Mat.
  Ital. A (7) \textbf{10} (1996), no.~1, 127--141. \MR{97a:54013}

\bibitem{CN}
Eduard {\v{C}}ech and Josef Nov{\'a}k, \emph{On regular and combinatorial
  imbedding}, \v Casopis P\v est. Mat. Fys. \textbf{72} (1947), 7--16.
  \MR{9,98e}

\bibitem{C}
M.~M. {\v{C}}oban, \emph{Operations over sets}, Sibirsk. Mat. \v Z. \textbf{16}
  (1975), no.~6, 1332--1351, 1372. \MR{54 \#7266}

\bibitem{D2}
Georgi~D. Dimov, \emph{On the commutability between hyperspaces and subspaces,
  and on {S}chmidt's conjecture}, Rend. Istit. Mat. Univ. Trieste \textbf{25}
  (1993), no.~1-2, 175--194 (1994), Proceedings of the Eleventh International
  Conference of Topology (Trieste, 1993). \MR{96d:54011}

\bibitem{D1}
\bysame, \emph{Some remarks on the commutability between hyperspaces and
  subspaces, and on {S}chmidt's conjecture}, C. R. Acad. Bulgare Sci.
  \textbf{47} (1994), no.~4, 5--8. \MR{1 332 595}

\bibitem{DV}
Georgi~D. Dimov and Dimiter Vakarelov, \emph{On {S}cott consequence systems},
  Fund. Inform. \textbf{33} (1998), no.~1, 43--70. \MR{2000b:03226}

\bibitem{E}
Ryszard Engelking, \emph{General topology}, PWN---Polish Scientific Publishers,
  Warsaw, 1977, Translated from the Polish by the author, Monografie
  Matematyczne, Tom 60. [Mathematical Monographs, Vol. 60]. \MR{58 \#18316b}

\bibitem{F}
Oskar Feichtinger, \emph{Hyperspaces with the kappa topology}, Topology Proc.
  \textbf{3} (1978), no.~1, 73--78 (1979), Proceedings of the 1978 Topology
  Conference (Univ. Oklahoma, Norman, Okla., 1978), I. \MR{80k:54011}

\bibitem{J}
Peter~T. Johnstone, \emph{Stone spaces}, Cambridge University Press, Cambridge,
  1982. \MR{85f:54002}

\bibitem{K}
Borislav Karaivanov, \emph{On the commutability between hyperspaces and
  subspaces}, Questions Answers Gen. Topology \textbf{14} (1996), no.~1,
  85--102. \MR{97a:54014}

\bibitem{M}
Ernest Michael, \emph{Topologies on spaces of subsets}, Trans. Amer. Math. Soc.
  \textbf{71} (1951), 152--182. \MR{13,54f}

\bibitem{S}
Hans-J{\"u}rgen Schmidt, \emph{Hyperspaces of quotient and subspaces. {I}.
  {H}ausdorff topological spaces}, Math. Nachr. \textbf{104} (1981), 271--280.
  \MR{84e:54011a}

\bibitem{Se}
Milan Sekanina, \emph{Topologies on systems of subsets}, General topology and
  its relations to modern analysis and algebra, IV (Proc. Fourth Prague
  Topological Sympos., Prague, 1976), Part B, Soc. Czechoslovak Mathematicians
  and Physicists, Prague, 1977, pp.~420--424. \MR{57 \#13826}

\end{thebibliography}
\providecommand{\bysame}{\leavevmode\hbox to3em{\hrulefill}\thinspace}
\providecommand{\MR}{\relax\ifhmode\unskip\space\fi MR }
\providecommand{\MRhref}[2]{%
  \href{http://www.ams.org/mathscinet-getitem?mr=#1}{#2}
}
\providecommand{\href}[2]{#2}

\end{document}